%% file: ghi_rever.tex
\documentclass[a4paper,11pt]{article}

\title{On certain generalized Hardy's inequalities and applications}

\author{Demetrios A. Pliakis}
\begin{document}

\maketitle

\newcommand{\mpalla}[1]{\ensuremath{B_{#1}}}
\newcommand{\mpallat}[1]{\ensuremath{B^*_{#1}}}
\newcommand{\eps}{\ensuremath{\epsilon}}
\newcommand{\mpallab}[1]{\ensuremath{\widehat{B}_{#1}}}
\newcommand{\kl}{\ensuremath{\nabla}}
\newcommand{\olwn}{\ensuremath{\int_{{\bf R}^n}}}
\newcommand{\ced}{\ensuremath{C_{\epsilon,\delta}}}
\newcommand{\sed}{\ensuremath{S_{\epsilon,\delta}}}
\newcommand{\sedj}{\ensuremath{S_{\epsilon,\delta,j}}}
\newcommand{\cehj}{\ensuremath{S_{\epsilon,\eta,j}}}
\newcommand{\red}{\ensuremath{R_{\epsilon,\delta}}}
\newcommand{\xw}[1]{\ensuremath{{\bf R}^{#1}}}
\newcommand{\pxw}[1]{\ensuremath{{\bf P}^{#1}}}
\newcommand{\kw}[1]{\ensuremath{C_{#1}}}
\newcommand{\ankw}[1]{\ensuremath{\widehat{C}_{#1}}}
\newcommand{\fo}{\ensuremath{\mbox{supp}}}
\newcommand{\pal}[1]{\ensuremath{\partial_\lambda^{#1}}}
\newcommand{\resoul}[1]{\ensuremath{R_0(\lambda)^{#1}}} 
\newtheorem{genhardy}[section]{Theorem}
\newtheorem{desalg}[section]{Theorem}
\newtheorem{desalglem}[section]{Lemma}
\newtheorem{lojas}[section]{Theorem}
\newtheorem{inhom}[section]{Proposition}
\newtheorem{essesa}[section]{Proposition}
\newtheorem{opes}[section]{Proposition}
\newtheorem{opesII}[section]{Proposition}
\newtheorem{salmv}[section]{Proposition}
\newtheorem{abg}[section]{Theorem}

\input{ghi_intro.tex}

\input{ghi_prel.tex}

\input{ghi_proof1.tex}

\input{ghi_proof2.tex}

\input{ghi_appl_1.tex}

\input{ghi_appl_2.tex}

\input{ghi_ref.tex}
\begin{center}
Author's address: TEI CRETE, Dept Electronics, \\
Chania, Crete, Greece\\
 email address: \ttfamily{dpliakis@gmail.gr}
\end{center}
\end {document}

%% file: ghi_intro.tex
\paragraph{Introduction}
The classical inequality of Hardy for smooth functions 
\( f  \in C_0^\infty(f\in\xw{}\setminus \{0\})\):
\[ 
\int_{{\bf R}} x^sf^2 \leq \frac{4}{(s+1)^2} \int_{{\bf R}} x^{s+2}(f')^2
 \]
for \( s \neq -1 \) can be generalized in various ways and provides a weighted
version of Poincar\'e's  inequality. The standard generalizations replace the 
weight  \( x^s\) with the radial variable or the boundary defining 
function of a smooth domain, are reduced  to  the one-dimensional 
case and are proved directly by partial integration, as far as the weight stays smooth.
 Here we we replace the weight by a homogeneous 
polynomial that is singular also away from the origin, 
so its zero set is a singular algebraic cone. In this case
 no direct method of the preceding form is available: the rectilinearization 
of such a set being non-trivial along its singularities. Specifically,
 singular algebraic varities are  
rectilinearized under the process of ''resolution of singularities'' then,
their singularities unfold and appear as ''normal crossings''. We follow this
procedure to the extent of '' reduction of multiplicity '' of an algebraic set 
and prove following generalization. 

Let $P(x_1,\dots,x_n)$ be a homogeneous polynomial of degree $d$ in
$n$-real variables belonging to the class ${\cal P}^{gH}$ that we define
 in the next paragraph.  Let $V(P)=\{ x\in {\bf R} ^n/ P(x)=0\}$ be the  
algebraic set that it defines. We introduce the Hardy factors:
\[
 {\cal H}^1(P)=P^{-\frac2d},\,\,\,\, {\cal H}^2(P)=\left|\frac{\nabla P}{P}\right|^2
\]
We prove the   following generalized Hardy inequalities \(\mbox{GHI}_i\):
\[
 \int_{{\bf R}^n}{\cal H}^i(P)f^2\leq C_i(P)\int_{{\bf R}^n}|\nabla f|^2
\]
for functions 
$f\in C_{0}^{\infty}({\bf R} ^n\setminus V(P))$. This inequality while it is 
elementary to prove when the algebraic variety $V(P)$ is smooth 
away from the origin, it is rather cumbersome when the variety is singular. 
The above inequality may be viewed as direct generalization of Hardy's 
 Here, we will consider 
the stratification of the algebraic variety
$V(P)$ by multiplicity and the inequality will be examined through the resolution 
of singularities process. This provides a finite covering, 
in every chart of which the algebraic set is reduced 
to normal crossings. The inequality is readily reduced to a 
corresponding one for inhomogeneous polynomials.

\begin{genhardy}
 
 Let \(P \)  be a nonhomogeneous polynomial of degree $d$. Then there is a  
constant \(c'_1(P)>0,i=1,2\) such that
\[
\int_{{\bf R}^n}{\cal H}^1(P) f^2 \leq c_1'(P) \sum_{i=1}^{h_1}||\nabla^i f||^2_2 
\]
If \( P\in {\cal P}^{gH}\) then we have the following:
\[
 \int_{{\bf R}^n}{\cal H}^2(P) f^2 \leq c_1'(P) 
\sum_{i=1}^{h_2}||\nabla^i f||^2_2
\]
\end{genhardy}

However, it is worthnoticing that we can refine the crude form of the
preceding inequality for an inhomogeneous polynomial of degree $d$, 
belonging in the class described below; precisely there is a constant
\(c_3(P)>0\) 
such that in the operator sense on the domain 
$C^{\infty}_0({\bf R} ^n\setminus V(P))$ there holds:
\[
 \int_{{\bf R}^n}{\cal H}^1(P) f^2 \leq c_3(P) \int_{{\bf R}^n}|\nabla f|^2+
(1+|x|^2)f^2 
\]

We present two applications of this inequality: 

\begin{itemize}

\item We apply the inequality to a particular case of a problem that 
motivated the study of these inequalities: the existence 
of an asymptotic expansion in powers and logarithms of the distributional 
trace of the heat operator corresponding to 
\[
H_{c,\alpha}=-\Delta + \frac{c}{|P|^\alpha}
\]
for $c \in C_{0}^{\infty}({\bf R} ^n)$ and small $\alpha >0$. 
The inequalities provide the required estimates for the domain, closure 
and the Neumann series of a suitable power of the 
resolvent of $H_{c,\alpha}$. The existence of the asymptotic expansion 
follows, in view of the singular asymptotics lemma \cite{c1},  
from the well known theorem of Bernshtein-Gelfand, \cite{bg}, 
for the meromorphic extension of integrals containing complex powers of 
polynomials.

\item We consider a smooth domain that approaches arbitrarily close a 
non-smooth one: it is defined by the level sets of a polynomial. 
In this domain  we compare the growth of functions to the growth of the 
polynomial defining it. This allows us to compare values of functions in the
following manner:   
in  a smooth metric in a Euclidean domain represent the domains of
given curvature growth by the semialgebraic sets defined by specific 
polynomials. Then we obtain estimates of the local growth of laplacian 
eignefunctions in terms of  the curvature growth. Here we prove the  
simple inequality that provides the 
growth of integrals of functions is such semilagebraic domains.

\end{itemize}
 
The article begins with a review of the local reduction of an algebraic 
set to normal crossings according to \cite{bm1} with certain comments that 
adjust the process with our purposes.  In the sequel we prove the preceding
theorem and its further generalizations and then we proceed to the 
applications.

%% file: ghi_prel.tex
\section{Local reduction of an algebraic set to normal crossings.}

We commence by reviewing the necessary definitions and results on blowing up
and local desingularization of an algebraic set with guide essentially
the presentation in \cite{bm1}.

\subsection{ Standard Constructions}

\subsubsection{Blown up space} In the various steps we will  use conical partitions of
 unity covering the euclidean ball centered at the origin
\[
 \mpalla{0,\epsilon,n}=\{x \in  \xw{n} / |x|^2:=x_1^2+\cdots x_n^2 < \epsilon\}
\] 
subordinate to the cones 
\[
\kw{\alpha;j}=\{x\in \xw{n} / x_j^2 > \frac{1}{1+\alpha}|x|^2\}
 \]
for \(\alpha >1\). Let \(\pxw{n}\) denote the 
$n$-dimensional projective space of lines through the origin in 
$\xw{n+1}$. Let 
\(\mpallat{0,\eps,n+1}=\mpalla{0,\eps,n+1}\setminus\{0\}\) be the punctured 
euclidean ball around the origin in $\xw{n+1}$, set
\[
\mpallab{0,\eps,n+1}=\overline{ \{ (x,l) \in
\mpallat{0,\eps,n+1} \times {\bf P}^n: x \in l \}}
\]
and let $\sigma : \mpallab{0,\eps,n+1} \rightarrow\mpalla{0,\eps,n+1} 
\,\,\,\, (x,l) \mapsto x$. Then $\sigma$ is proper, restricts to a homeomorphism 
over $\mpallat{0,\eps,n+1}$, and $\sigma^{-1}(0)={\bf P} ^n$. This mapping is 
called the blowing up of $\mpalla{0,\eps,n}$ with center $\{0\}$. 
In a natural way, $\mpallat{0,\eps,n+1}$ is an algebraic submanifold
of $\mpalla{0,\eps,n} \times \pxw{n}$:

\paragraph{Coordinates} Let \((x_1,\dots,x_{n+1})\) denote the affine coordinates in 
\(\xw{n+1}\) and let \(t = [t_1:...:t_{n+1}]\) denote the homogeneous 
coordinates of $\xw{n}$. Then 
\[\mpallab{0,\eps,n+1}=
\{ (x,t) \in \mpalla{0,\eps,n+1} \times \pxw{n}: x \wedge t =0\}\]
Furthermore $\mpallab{0,\eps,n}$ is covered for $\alpha>1 $
by the conical charts \(j=1,\dots,n+1\):
\[
\ankw{\alpha:j} =
\{ (x,t) \in \mpallab{0,\eps,n+1}: t_i^2 > \frac{1}{1+\alpha}|t|^2\}, 
\]
 with coordinates $(x_{1,i},...,x_{n+1,i})$, for each $i$, where 
$$x_{ii}=x_i, \hspace{1truecm} x_{ji}=\frac{t_j}{t_i}, i \neq j$$
with respect to these local coordinates, $\sigma$ is given by
$$x_i=x_{ii} \hspace{1truecm} x_j=x_{ii}x_{ji}, \hspace{0.2truecm} i \neq j$$
Let $n>c$ and $B_\epsilon^{n-c}(0)$ then the
mapping $\sigma \times id : \mpallab{0,\eps,c} \times 
\mpalla{0,\eps,c}
\rightarrow B_\epsilon^c(0) \times B_\epsilon^{n-c}(0) $ is called the 
blowing up of $B_\epsilon^c(0) \times B_\epsilon^{n-c}(0)$ 
with center $C:=\{0\} \times B_\epsilon^{n-c}(0) \subset {\bf R} ^n$
and it is denoted by $\mbox{Bl}_{C}(B_\epsilon^c(0)\times B_\epsilon^{n-c}(0))=
\mpallab{0,\eps,c}\times B_\epsilon^{n-c}$.

\subsubsection{Blown up volumes and vector fields} 
Let the usual volume in $\xw{n}$ be denoted by 
 $v_n$, then under blow up with center of codimension $c>1$  
considered in the $i$-th chart it pulls back the volume 
\[
v_n=x_i^{c-1}\widehat{v_n} 
\]
It is noteworthy the way that the vector fields that generate 
dilation transform under blowing up or down. Let then $x=\sigma(y)$ and
\begin{eqnarray*} 
 D_{x_i}=x_i\frac{\partial}{\partial x_i},\,\,\,\, 
D_{y_i}=y_i\frac{\partial}{\partial y_i},\\ 
E_{0}=\sum_i x_i\frac{\partial}{\partial x_i},\,\,\,\,\, 
E_{1}=\sum_i y_i\frac{\partial}{\partial y_i}
\end{eqnarray*}
then in the $k$-th chart we have the formulas:
\begin{eqnarray*}
i\neq k:\,\,\,D_{x_i}=D_{y_i} , \\ D_{k_1} =E_0, \\
 E_1=2E_0-D_{x_k}
\end{eqnarray*}
Notice theat the Euler vector field \(E\)  
is expressed in the radial variable \(r=|x|\) as 
\[
E=r\partial_r=D_r
\]

We will consider mappings obtained as a finite
 sequence of local blow-ups ; i.e., 
$\pi_N = \sigma_1 \circ ...\circ \sigma_N$,         
where for each  $i=1,...,r,  
\sigma_i : \widehat{V}_i \rightarrow V_i$ is a local
blow up with center $C_i=\{0\} \times \mpalla{0,\eps,n-c}$
of the preceding form with  
$\widehat{V}_i=\mpallab{0,\eps,c_i}\times \mpalla{0,\eps,n-c_i}$ and 
$V_i=\mpalla{0,\eps,c_i}\times \mpalla{0,\eps,n-c_i}$.

\paragraph{The conical atlas partition of unity} 
We conclude with the partition of unity of the 
punctured ball \(\mpallat{0,\eps,n+1}\) subordinate to its conical covering.
Let \(\varphi\in C_0^\infty(\xw{}_+)\) with \(\fo(\varphi)\subset 
[0,1+\varepsilon), \phi\equiv 1 \) in \([0,1]\) then set
\[
\chi_j(x) = \varphi\left(\frac{r}{\sqrt{1+\alpha}|x_j|}\right)
\]
We compute its derivatives:
\[
|\kl^\ell \chi_j|\leq \frac{1}{r^\ell}\sum_{i_1+\cdots+i_\ell=\ell}
|\varphi^{i_1}|\cdots|\varphi^{i_l}|\leq \frac{C_\ell}{r^\ell}
\]
This formula is very important because in the end of the blow-up
process we encounter the derivatives of the localizations functions.
These will require the further application of the one dimensional Hardy
inequality, i.e. in the radial variable.

\subsection{The local desingularization algorithm.} Here we 'll  follow the
proof of the local desingularization theorem in algorithm devised in [BM]
and developed in the conventions that we need for the inequalities. 
It consists of two steps the determination of the center and 
the reduction of the multiplicity.

\vspace{0.5truecm}

\begin{desalg}
 
 Let \(P : \xw{n} \rightarrow \xw{}\) be a regular 
function. Then there is a countable collection of regular mappings 
$\pi_r: W_r \rightarrow {\bf R} ^n$ such that:

\begin{enumerate}

\item  Each $\pi_r$ is the composition of a finite 
sequence of local blows up (with smooth centers)

\item There is a locally finite covering $U_r$ of ${\bf R} ^n$ such that
\(\pi_r(W_r) \subset U_r\) for all $r$.

\item If $K$ is a compact subset of ${\bf R} ^n$, there are compact
subsets $M_r \subset W_r$ such that $K= \bigcup_r \pi_r(M_r)$.
The union is finite by (2).

\item  For each $r, P \circ \pi_r$ is locally normal crossings.

\end{enumerate}

\end{desalg}

\paragraph{Determination of the center.} 
Let $a \in V(P),mbox{ord}_a(P)=:m$ and choose coordinates $x=(x_1,\dots,x_n)$ such 
that $x(a)=0$ and, $\mbox{in}_a(P)(0,...,0,x_n) \neq 0 $ the lowest degree 
homogeneous component of the polynomial. 
Moreover denote by ${\tilde x}=(x_1,...,x_{n-1})$ and let $h_j=V(x_j)$ be 
the coordinate hyperplanes.  Since $(\partial^m_{x_n}P)(x) \neq 0 $ then  
\[
(\partial^{m-1}_{x_n}P)(x) \sim x_n - H({\tilde x}) =: x_n'
\]
for some regular function $H$. Then we perform the division 
\[
P(x)=Q(x)x_n^m + \sum_{0 \leq k < m} c_k({\tilde x})x_n^k
\] 
Moreover, possibly after translation, we may assume that 
$c_{m-1}({\tilde x}) = 0$ and also observe that each 
$c_k({\tilde x}) = \partial^k_{x_n}|_{h_n}P, 
{\tilde x}_n \sim \partial_{x_n}^{m-1} P $
  
\vspace{0.25truecm}
\noindent Let  $ a \in h_n$ then for $0 \leq k < \mbox{ord}_aP= \mu _P(a)$
introduce the sets 
\[
  {\cal P}(a):= \{ (P,\mbox{ord}_aP)\}, \hspace{0.2truecm}
  {\cal C}_P(a) := \{ (c_{k}, \mu_P-k)\}
\] 
The union of stata of multiplicity  at least $\mu_P(a)(=m) $  is denoted by 
\[
S_{{\cal P}(a)}:= \{ x : \mbox{ord}_x P \geq \mu_P(a) \}
\] 
as well as that 
\[
S_{{\cal C}_P(a)}:=\{ x: \mbox{ord}_xg \geq p,\,\,\,
  \mbox{forall}\,\,\, g \in {\cal C}_P(\alpha)\}
\]
First we will use an induction on $n$ (and on $m$), 
to arrive at the particular instance when for all $k$
\[
(*)\,\,\,\,\,\,\,\,c_k({\tilde x})=({\tilde x}^{\gamma})^{m-k} 
c_k^*({\tilde x})
\]
while  $\gamma \in \frac{1}{m!} {\bf N} ^{n-1}$ and for some 
$k_0, c_{k_0}^*(0) \neq 0. $ In order to handle at once the various 
$c_k$'s we define the auxiliary function 
\[
A_P({\tilde x}) := \mbox{ product of all non zero}\,\,\,\,
c_{k}^{\frac{m!}{m-k}}\,\,\,\, \mbox{and all their nonzero differences}
\]
The inductive assumption asserts that there is already a 
uniformization for \(A_P:\) 
\[
A_P \sim x_1^{a_1}...x_{n-1}^{a_{n-1}}
\]
This implies in first place that each nonzero 
$ c_k({\tilde x})=({\tilde x}^{\Omega_k}) c_k^*({\tilde x})$
with $\Omega_k \in  {\bf N} ^{n-1}$ and $c_k(0)^* \neq 0 $. Moreover each 
nonzero 
$c_k^{\frac{m!}{m-k}}-c_j^{\frac{m!}{m-j}} \sim {\tilde x}^{\Lambda_{ij}}$, 
with ${\Lambda}_{ij} \in     ({\bf N} ^{n-1})^*$. 
The following elementary  lemma suggests:

\vspace{0.5truecm}

\begin{desalglem}

 Let $x=(x_1,...,x_n)$.
If $a(x)x^{\alpha}-b(x)x^{\beta} = c(x)x^{\gamma}$ and 
$a(0)b(0)c(0) \neq  0$ then either $\alpha \in \beta + {\bf N} ^n$ or
$\beta \in \alpha + {\bf N} ^n$

\end{desalglem}

\noindent The lemma implies that the set 
${\cal E} := \{ \frac{1}{m-k} \Omega_{k} \} \subset \frac{1}{m!} {\bf N} ^{n-1}$ 
is totally ordered with the induced partial ordering from ${\bf N} ^{n-1}$ 
and therefore there exist a $\rho = \min({\cal E}) \in \frac{1}{m!} {\bf N} ^{n-1}$.
 Therefore we are reduced to the case (*) with    
$\gamma = \rho \in  \frac{1}{m!} {\bf N} ^{n-1}$. 
We show that the special case (*) implies reduction in multiplicity by 
blowing up succesively the components of $S_{\cal P}(a)$.
\[
S_{{\cal P}(a)}= S_{{\cal C}_{P}(a)}=
 \{ x : x_n=0 , \,\,\,\ \mbox{ord}_{x} ({\tilde x}^{\gamma}) \geq 1 \}
= \bigcup _I Z_I
\]
where $Z_I=\bigcap_{i\in I}(h_i\cap h_n)$ and 
 $I \subset \{1,...,n-1\}, \mbox{card}(I)=\nu_m-1 $ minimal such that 
$\sum_{i \in I} \gamma_{i} \geq 1$ or equivalently that 
$0 \leq \sum_{i \in I} \gamma_{i} -1 < \gamma_{k}$ for all $k \in I$.
Actually, these serve as centers of the desingularizing blowups
$C= Z_I$,  for all $k$ we have that $$ \mbox{ord}_{C}c_{k} \geq m-k > 0, $$
where $\mbox{ord}_{C}c_{k} = \inf_{x \in C}\bigl( \mbox{ord}_x c_k \bigr)$.

\vspace{0.2truecm}
\paragraph{Reduction of the multiplicity.}
Let $a \in V(P)$ and choose coordinates such that $x(a)=0,I=(1,\dots,\nu_m-1)$. Let further
\begin{eqnarray*}
\kw{\alpha;\epsilon;j,n,\nu_m}=
\left(\kw{\alpha;j,n} \cap \mpalla{0,\eps,\nu_m}\right) \times \mpalla{0,\eps,n-\nu_m}\\
\ankw{\alpha;\epsilon;n,\nu_m}=
\left(\ankw{\alpha;\epsilon;j,n,\nu_m} \cap \mpalla{0,\eps,\nu_m}\right) \times 
\mpalla{0,\eps,n-\nu_m} 
\end{eqnarray*} 
and the blowup 
\[ 
\sigma : \mpallab{0,\eps,n-\nu_m}\times\mpalla{0,\eps,n-\nu_m} 
\rightarrow \mpallab{0,\eps,\nu_m}\times \mpalla{0,\eps,n-\nu_m}
\]
at $\beta \in \sigma^{-1}(0)$. We calculate in these conical charts:

\begin{itemize}
\item  Let $\beta \in \ankw{\alpha;\epsilon;n,n,\nu_m}$ 
then the  strict transform $\sigma|C_n $ maps as follows
$$ x_n=y_n,  \hspace{0.2truecm} x_j=y_jy_n, j=1,\dots,\nu_m-1, \hspace{0.2truecm}
x_s=y_s, n>s>\nu_m.$$  Its effect on the polynomial is 
 \[ (\sigma^*(\chi_nP))(y)  =  y_n^mP^1_n(y)=y_n^m
[(\sigma^*Q)(y)+ 
\sum_{0 \leq k < m} (\sigma^*c_k)(y)y_n^{k-m}]
\]
and we observe that $Q(\sigma(\beta)) \neq 0$ while 
$\bigl(\sum_{0 \leq k < m}c_{k}(y)y_{n}^{k-m}\bigr)(\sigma(\beta))=0$
and hence $P_n^1(\beta) \neq 0$.

\item In the conical sector 
\(\ankw{\alpha;\epsilon;j,n,\nu_m},\,\,\,j \in I\)
the strict transform $\sigma$ maps as follows  
\[ x_n=y_n, \hspace{0.2truecm} y_j,x_j=y_j, \,\,\,\,\, x_k=y_ky_j,
j\neq k=1,\dots,c-1, \,\,\,\,\, x_s=y_s, s>n.
\]
Then the polynomial becomes 
\[
(\sigma^*(\chi_jP))(y) = y_j^mP^1_j(y)=y_j^m\bigl(
(\sigma^*Q(y))y_n^m+ 
\sum_{0 \leq k < m} (\sigma^*c_k )({\tilde y})y_j^{k-m}
y_n^k \bigr) 
\]
since in this sector
\[
{\tilde x}= \sigma ({\tilde y}), {\tilde y}=(y_1,...,y_{n-1}).
 \]
Observe that $\partial^{m-1}_{y_n} P^1_j \sim y_n$
since $Q(\alpha) \neq 0$ while  $h_n'=\sigma^{-1}(h_n)=\{y_n=0 \}$ and 
$ c_{m-1} = 0$ identically. 

\end{itemize}

We conclude that for all points on the exceptional 
divisor the order is not bigger than $m$, $\beta \in \sigma^{-1}(C) , 
\mbox{ord}_{\beta}P^1
 \leq d$ if  $P^1$ denotes the resulting regular function.
Therefore assume that $\mbox{ord}_{\beta}P^1 = m$ iff 
$\beta \in h_n \bigcap \{ \mbox{ord}_{y}C'_{k} \geq m-k \}$ where 
$$c'_k=y_j^{k-m}(\sigma^*c_k)({\tilde y})= 
({\tilde y}^{\gamma'})^{m-k}.(c_k^*(\sigma({\tilde y}))$$
and ${\tilde y}^{\gamma'}:=y_j^{-1}.({\tilde x}^{\gamma}\circ \sigma)$ and 
there is $k_0, (c_{k_0}^*)(\sigma(\beta)) \neq 0$
by the particular case (*).
Hence, $\gamma_i'=\gamma_i$, if $i \neq j $ and $\gamma_j'=
\sum_{i \in I}\gamma_i -1 $ therefore $1 \leq |\gamma'| < |\gamma|$
and since $|\gamma|,|\gamma'| \in \frac{1}{m!} {\bf N} ^{n-1}$
it follows that after no more
than $|\gamma|m!$ blows up of this type multiplicity has to decrease.

\paragraph{Remarks.} 

\begin{itemize}

\item Let $P \in {\bf R} [x_1,\dots,x_n]$ and $V(P):=V$ the variety it 
defines. Set then 
$$P^k(x)=\sum_{j_1+\dots+j_n=k}\bigr(\frac{\partial^kP}{\partial x_1^{j_1}
\dots \partial x_n^{j_n}}(x)\bigl)^2$$ and consequently $V_k:=V(P_k)$
The stratification by multiplicity we 'll use consists of the strata that
are semialgebraic sets:
$$ \Sigma_k=(V_k\setminus V_{k+1})\cap V$$ 
If $P$ is homogeneous then by Euler's theorem 
$\Sigma_{k+1} \subset \overline{\Sigma_k}$ while the implicit function 
theorem asserts that $\Sigma_k$ are smooth.

\item Let $I(V)$ be the sheaf of germs of regular functions on ${\bf R}^n$ 
that vanish on $V$. Let $x \in \Sigma_k$ then set
$$L_{x,k}= \{ \xi \in {\bf R} ^n :  P\in I_x(V),
\sum_j\xi_j\frac{\partial P}{\partial x_j} \in I_x(V) \}$$ and also that
$E_{x,k}=L^{\perp}_{x,k}$. 

\item  The class ${\cal P}^{gH}$ of polynomials consists of those $P$ such 
that when  $codim\Sigma_k=2$ for some $k$ then 
$$\{(x+E_{x,k})\cap V_x \} \setminus \{ x\} \neq \emptyset $$ where $V_x$
is the germ of $V$ at $x$.

\item If $P \geq 0 $ everywhere as well as that $P_1,\dots,P_p \geq 0$ 
everywhere and it is true that
$$P=\sum_j P_j $$ then $V(P)$ is not a hypersurface. Then the same
procedure with more complicated details brings the set to normal crossings, cf 
[BM]. However, if we assume that each of the $P_j$'s belongs to the class 
${\cal H}$ then  we can proceed without appealing to the desingularization 
for codimension$>1$.

\item We would assume that \(P\) is irreducible otherwise 
split it in its factors and use Young's inequality to deal with its factor
separately.

\end{itemize}

\paragraph{The inequalities of Lojasiewicz.} Through the paper certain 
distances from algebraic sets are estimated by the values of the defining 
polynomials through the fundamental Lojasiewicz inequalities, cf [BM]:

\begin{lojas}
 
Let $P$ be a regular function on an open subspace 
$M \subset {\bf R}^n$. Suppose that $K$ is a compact subset of $M$, 
on which $V(|\nabla P|^2) \subset V(P)$. Then there exist $c,c'>0$ and $\mu, 
0 < \mu \leq 1,\nu>1$ such that 
\[
|\nabla P(x)| \geq c|f(x)|^{1-\mu},\,\,\,\, |P(x)| \geq c' d(x,V(P))^{\nu}
\]   
in a neighborhood of $K, \sup(\mu),\sup(\nu) \in \bf Q$.

\end{lojas}

\noindent 
It is true that  $ \nu(K) \leq \sup_{a\in K}\mu_a(P)$. 
Furthermore in the case of a homogeneous polynomial function in
${\bf R} ^n$ the constant depends on conical neghborhoods of the 
origin.

In the case of a homogeneous polynomial \(P\) of degree \(m\) then
the algebraic cone decomposes as \(V(P)=\xw{}_+\times K\) where
\( K\) is the trace of the cone on the sphere. Then Lojasiewicz inequality 
suggsets that for \( C(K)>0,\mu=\frac1m\):
\[
 |\kl P|\geq C_m|P|^{1-\mu}
\]
This follows from the application of Lojasiewicz on the trace, under stereographic 
projection in coordinates \((r,\xi)\in \xw{}_+\times \xw{n-1}\), and for applying Young's 
inequality for \(p,q>1, \,\,pq=p+q,\):
\begin{eqnarray*}
 |\kl P(x)|\geq mr^{m-1}|\tilde{P}(\xi)|+r^m|\widehat{\kl P}(\xi)|\geq\\
\geq mr^{m-1}|P|+C r^m|P|^{1-\mu}\geq C'_mr^{\frac{1}{\mu p}-\frac{m}{q}}
|P(x)|^{1-\frac{p}{\mu}}
\end{eqnarray*}
We choose as \(p=\frac{m+\mu}{m}\) scale and arrive at the result \(\mu=\frac1m\).

%% file: ghi_proof1.tex
\section{The inequalities.}

Here we will prove the inequalities \( \mbox{GHI}_1,\mbox{GHI}_2\)  
for an inhomogeneous polynomial \(P\) of degree
\(d \) in the class \({\cal H}\) and an arbitrary function \(f\in 
C^\infty_0({\bf R}^n\setminus V(P))\):
\[
 I^i[P](f)=\int_{{\bf R}^n} {\cal H}^{i}(P) f^2 \leq C \int_{{\bf R}^n} 
\sum_{j=0}^{h_i} |\nabla^jf|^2=||f||_{H^{h_i}({\bf R}^)}, \,\,\,\,\,\,\,\,\,
(\mbox{GHI}_i)
\]
where 
\[
 {\cal H}^1(P)=P^{-\frac2d},\,\,\,\, {\cal H}^2(P)=\left|\frac{\nabla P}{P}\right|^2
\]
This inequality will be based on the fact that after a suitable number of blow-ups 
with suitably chosen centers the multiplicity of the polynomial has to decrease. 
The choice of centers is provided in the proof of the local desingularization theorem. 
 
To fix the ideas assume that we are localised in a tubular neighbourhood of \(V(P)\)
of  width $\epsilon^{\frac1d}$:
\[
N_{\epsilon}(P)=\{ x\in {\bf R} ^n/ |P(x)|<\epsilon\}
 \]
and also the tubes of width  \(\epsilon_k^{\frac1d}\)  that enclose
the strata of multiplicity \(\Sigma_k\):
\[
N_{\epsilon_k,k}(P)=\{ x \in {\bf R} ^n/Q_k(x)=\sum_{j=0}^k P^j(x)<\epsilon_k^2\}
\]
We assume that \(\epsilon_1,\dots,\epsilon_m,\epsilon \) 
are chosen so that this system of tubes \(N_1,\dots,N_m\) exhaust \(N_\epsilon\)
and using functions of the form \(\chi\left(\frac{Q_k}{\epsilon_k}\right)\),
for \(\chi\) a one-dimensional cut-offs we localize in these sets. Thanks to
Lojasiewicz these cut-off functions when differentiated stay away from the variety.
The integral then is splitted up as
\[
 I^i[P](f)=\sum_{j=0}^m I^i_j[P](f)
\]
where 
\[
 I_j^i[P](f)=\int_{N_j} {\cal H}^i(P) f^2
\]
As a matter of fact we have that 
\[
 \fo(f)= \bigcup_{j=0}^m (\mbox{supp}(f)\cap N_j)
\]
We are ready to set-up the multiplicity reducing algorithm. Select then points 
\[
 a \in V(P)\cap N_{m,\epsilon}(P):\,\,\,  \mbox{ord}_aP=m,\,\,\,\,
\mbox{codim}(S_{P(a)})= \nu_m
\] 
and choose a system of coordinates such that $x(a)=0$ and also that 
\[
P(0,\dots,0,x_n)\sim x_n^m 
\] 
The initial change of coordnates that rectifies the center we will
denote by \( K_m:N_m \rightarrow {\bf R}^{n-\nu_m}\times {\bf R}^{\nu_m} \),
which consists  of algebraic maps localized through \(\varphi_{m,\ell_m}\)
at such  points on \(a\in S_m\) :
\[
 K_m=\sum_{\ell_m=1}^{N_m} \varphi_{m;\ell_m} K_{m;\ell_m}
\]
These maps have jacobians \(0<c_m(\epsilon_m)\leq J(K_{m;\ell_m})\leq C_m(\epsilon_m)\)

\subsection{The inequality \(\mbox{GHI}_1\) for \({\cal H}^1=P^{-\frac2d}\)}

\paragraph{ The inequality in the multiplicity reduction step \(m \Rightarrow m-1\)} 
Using the conical partition of unity \(\{ \chi_{m,k}\}\) subordinate to the covering
 in
\[
 \bigcup_{k=1}^{\nu_m} K_m(N_{\epsilon,m}(P)) \cap {\bf R}^{n-\nu_m}\times C_k
\]
we localize and compute that  

\begin{equation}
I_m(P)[f]=\int_{N_{\epsilon_m,m}} {\cal H}^1(P)f^2(x) v_n
\leq c\sum_{k=1}^{\nu_m} 
\int_{N_{m;1,k}}
\left|D_k(|P_{1,m,k} |^{-\frac 1d}f_{m;1,k})\right|^{2}
\left({\cal D}_k^1\right)^2
\hat{v}_n 
\end{equation} 
where
\begin{eqnarray*}
c= \kappa_m(\nu_m,d)=\left(\frac{16d^2}{d(\nu_m+1)+m)}\right)^2\\
N_{m;1,k}=\sigma^{-1}(N_{m;0,k}),\,\,N_{m,0,k}=K_m(N_{\epsilon_m})\cap 
{\bf R}^{n-\nu_m}\times C_k\\
D_k=x_k\partial_{x_k}, \,\,{\cal D}_k^1=x_k^{-\frac{m}{d}+\frac{\nu_m-1}{2}}\\
f_{m,1,k}=\chi_{m,k}f
\end{eqnarray*}
The last integral is majorized further by:
\[
\int_{N^1_{\epsilon_m,m}}
 {\cal H}^1(P_{m;1,k})\left[|D_kf^1_k|^2 +\left(Q^m_{1,k,1}(P)
f_{m;1,k}\right)^2\right] 
\left({\cal D}_k^1\right)^2 \hat{v}_n 
\]
The singular term is made of the following factors that we are going to keep track 
in the process: 
\begin{eqnarray*}
Q^m_{1,k,1}(P)=D_k(\log|P^{m,1}_k|),\\
 {\cal H}^1_{m;1,k} = \left(P_{m;1,k}\right)^{-\frac2d}, \,\,\,\, 
P_{m;1,k} =x_k^{-m}\sigma^*(\chi_kP)
\end{eqnarray*}
We perform successive blowups of this type in order to  reduce the multiplicity of 
the algebraic set and every time  apply the one dimensional inequality. 
Therefore,  we introduce inductively for the blow up years \(j>1\) and 
corresponding blow - up chart \(k_j\) the following functions that can come out:
\begin{eqnarray*}
N_{m;j} = \sigma^{-1}(N_{m;j-1}), \,\,\,\,\, N_{m;0}=N_{\epsilon,m}\\
f_{m;0,k,0}= \chi_{m,k} f,\,\,\,\,\,
 f_{m;j,k_j,l}=\sigma^*(\chi_{k_j}f_{m;j-1,k_{j-1},l}), 
\\
f_{m;j,k_j,l}= D_{k_j}\left(f_{m;j,k_{j-1},l-1}\right)\\
{\cal H}^1_{m;j,k_j}= x_{m,k_j}^{\frac{2m}{d}}
\sigma^*(\chi_{m,k_j}{\cal H}^1_{m;j-1,k_{j-1}}),\\
 {\cal D}_{m;j,k_j,l}=D_{k_j}({\cal D}_{m;j,k_{j-1},l-1}),\,\,\,\\
{\cal D}_{m;j,k_j,l}=\sigma^*(\chi_{m,k_j}{\cal D}_{m;j-1,k_{j-1},l})\cdot
{\cal D}_{k_j}^1,\\
Q^m_{j,k_j,l}=D_{k_j}(Q^m_{j,k_j,l-1})\\
Q^m_{j,k_j,l}=\sigma^*(\chi_{m,k_j}Q^m_{j-1,k_{j-1},l})
\end{eqnarray*}
The integral is then majorised by the sum after \(\gamma_m\) generations of blow-ups
that are neccessary for the multiplicity to reduce to \(m-1\):
\[
I_m(P)[f]\leq C(\epsilon_m,m,d,\nu_m) \int_{N_{\gamma_m}}\sum_{k_{\gamma_m}=1}^c
{\cal H}^1_{m;\gamma_m,k_{\gamma_m}}
\Phi_{m,1;k_{\gamma_m}}^2
\] 
for the functions that encode the history of blow ups
\[
\left(\Phi_{m,1;k_{\gamma_m}}\right)^2=\chi_{k_{\gamma_m}}^2\left[\sum_{l_1+l_2+l_3=\gamma_m}
\left(
f_{m;\gamma_m,k_{\gamma_m},l_1}{\cal D}_{m;\gamma_m,k_{\gamma_m},l_2}
Q^m_{\gamma_m,k_{\gamma_m},l_3}\right)^2\right] 
\]
This sum extends over the set of all \(\nu_m\)-adic 
numbers  with  \(\gamma_m\) digits: \( \Lambda_{\gamma_m}(c)\).
After these blow ups the polynomial \( P_{m;\gamma_m,k_{\gamma_m}}\) 
has  multiplicity \(m-1\).  
In the next \(\gamma_{m-1}\) generations of blow-ups:
\[
 I_m(P)[f]\leq  c\sum_{k_{\gamma_m}=1}^{\nu_m}
I_{m-1}[P_{m;\gamma_m,k_{\gamma_m}}](\Phi_{m,1;k_{\gamma_m}})
\]
Therefore we have to work  with
\[
 I_m[P](f)+I_{m-1}(P)[f]\leq c\int_{N_{m;\gamma_m}}{\cal H}_1(P_{m;\gamma,k_{\gamma_m}})
\Phi_{m,1;k_{\gamma_m}}^2+
\int_{N_{\epsilon_{m-1},m-1}} {\cal H}_1(P) f^2
\]
and proceed analogously. The desingularization algorithm guarantees that 
on the set \(N_{m;\gamma_m}\cap N_{\epsilon_{m-1},m-1} \)  
the polynomials \( P_{m;\gamma_m,k_{\gamma_m}},P\) have
 multiplicity \(m-1\). Therefore the choice of blow up centers entails 
the change of coordinates \(K_{m-1}\) 
that we trace in the summands. Their appearance modifies the
constant \(c\). Thus we have \(m\) generations of \(\gamma_1,\dots,\gamma_m\) years
of blow-ups.

\paragraph{Summing up for \(m=1\) and final step}
We arrive at the \(m=1\) stratum with the following sum of integrals:
\[
 I_1(P)[f]=\sum_{1\leq i_2\leq i_1\leq m}\int_{N_{i_1,i_2}} 
{\cal H}^1(P_{i_1;i_2,k_{i_2}}) \Phi^2_{i_1,i_2;k_{i_2}} 
\]
where the pair of indices stands keeps track of the origin of the function 
\(\Phi_{i_1,i_2;k_{i_2}}\) while the polynomials \(P_{i_1;i_2,k_{i_2}}\) 
have multiplicity \(1\).  Hence we change coordinates by the map \(K_1\)
and conclue with an applcation of Hardy's inequality.
 
\paragraph{Return to \(N_\epsilon\)} In order to return back to \(N_\epsilon\).
We summarize the process that we followed:

{\it Rectilinearization of the centers, blow up till we reduce the multiplicity by 1,
again rectilinearization of the new centers and new blow ups etc}

This finally made up a map: \( {\cal B}:N_\epsilon\rightarrow \tilde{N}_\epsilon\),
which is a piece-wise algebraic diffeomorphism outside a variety of positive 
codimension originating from the exceptional divisors of the blow ups.  
Clearly the process can be reversed and the formula derived above. 
We will perturb back the centers and blow down. The rectilinearization maps
just modify the constants back. The terms that have been produced will be examined then.
We examine now the \(\Phi_{i_1,i_2}\) terms: the \(Q\) factors are constants hence 
they just modify the consntant coefficients.  The \({\cal D}\)-terms split in two
terms: those that consist of the jacobians and the ''Hardy divisor factor''
\(y_k^{-2\frac{i_1}{d}}\). In the course of the blow down process the jacobian terms
dissapear while in the conical charts of  the \( {\bf R}^{\nu_m}\) factor 
of the coordinate system we have that for the radial variable \(r\)  :
\[
 y_k^{-2\frac{i_1}{d}}\leq c r^{-2\frac{i_1}{d}}
\]
 which  combine with the blow down formula to give us that
\[
 |D_{y_k}(f\sigma)|^2= |E(f)|^2\leq r^2 |\nabla f|^2
\]
If \( r^{-\beta}\) persists then we apply again the usual Hardy inequality.
The blowing down process will effect the replacement of the 
\( \left|D\left(f{\cal D})\circ \sigma\right)\right|\) with a term  \( |\nabla f|\)
and finally we get the result stated in the introduction

\subsection{The inequality \(\mbox{GHI}_2\) for 
\({\cal H}^2(P)=\left|\frac{\nabla P}{P}\right|^2\).} 
Here we separate at each year of blow up the  ''Hardy divisor term'', \(x_k^{-2}\)
 and blow down directly, which could have been done also in the  preceding case.
We apply the following elementary generalization of the Hardy's inequality refered in 
the introduction, for \(\ell \in {\bf N}, f\in C_0^\infty({\bf R}\setminus \{0\})\):
\[
 \int_{{\bf R}}\frac{f^2}{x^{2\ell}}\leq C_\ell \int_{{\bf R}} 
\left(f^{(\ell)}\right)^2 ,\,\,\,\,\, \mbox{(EGHI)}
\] 
and the appropriate detemination of the factors \({\cal D}, Q\) that appear during the 
process.  Therefore we start with the first blow up in the conical charts 
and obtain in view of
\[
 {\cal H}^2(P)  \leq \frac{m^2}{x_k^2}+\frac{4}{x_k^4}+
4\left({\cal H}^2(P_{m;1,k}\right)^2, 
\]
that since \(m\geq 2\)
\begin{eqnarray*}
 I_m(P)[f]=\int_{N_{\epsilon_m,m}} {\cal H}^2(P)f^2(x) v_n
\leq c\sum_{k=1}^{\nu_m} 
\int_{N_{m;1,k}} \left(\left|D_kf_{m;1,k}\right|^{2}+
\left|D_k^2f_{m;1,k}\right|^{2}\right)({\cal B}_{m;1,k_1})^2\\
+c\sum_{k_1=1}^{\nu_m}
\int_{N_{m;1,k_1}}\left({\cal H}^2_{m;1,k_1}\right)^2f_{m;1,k_1}^2
({\cal B}_{m;1,k_1})^2
 \hat{v}_n 
\end{eqnarray*} 
where
\begin{eqnarray*}
c=\kappa_m(\nu_m,d)=\frac{16m^2}{(2\nu_m+1)^2}\\
N_{m;1,k_1}=\sigma^{-1}(N_{m;0,k_1}),\,\,N_{m,0,k_1}=K_m(N_{\epsilon_m})
\cap \left({\bf R}^{n-\nu_m}\times C_{k_1}\right)\\
f_{m;1,k_1}=\sigma^*(\chi_{m,k_1}f)\\
D_{k_1} f_{m;1,k_1}=\frac{1}{x_{k_1}}E(f), \,\,{\cal B}_{m;1,k_1}=
|x_{k_1}|^{\frac{\nu_m-1}{2}},\\
 {\cal H}^2_{m;1,k_1} = \left(\nabla \log|P_{m;1,k_1}|\right)^2, \,\,\,\, 
P_{m;1,k_1} =x_{k_1}^{-m}\sigma^*(\chi_{k_1}P)
\end{eqnarray*}
The first terms is  blown down directly since in the conical chart 
\[ 
C_k: |x_k|^2\geq  \frac{\alpha^2}{1+\alpha^2}r^2 
\] 
and hence implies that
\[
 |D_kf_{m,k,1}|\leq \frac{\alpha^2}{1+\alpha^2}|\nabla f|^2
\]
Therefore we have by an application of Hardy to the derivative of the jacobian
\({\cal B}_{m;1,k_1}\) we could comprise all terms in the following inequality
\[
 I_m[P](f)\leq C\left(||\nabla f||_{H_1(N_m)}^2+I_{1;m;j}[P](f)\right) 
\]
where we have the usual Sobolev space norm:
\[
||\nabla f||_{H^1(N_m)} =\int_{N_{m}}|\nabla f|^2+|\nabla^2f|^2 
\]
and 
\[
 I_{1;m;1}[P](f)=\sum_{k=1}^{\nu_m} \int_{N_{m;1,k}}
\left({\cal H}^2_{m;1,k}\right)^2
\left(f_{m;1,k}{\cal B}_{1,k}\right)^2\hat{v}_n 
\]
Set then
\begin{eqnarray*}
 f_{m;j,k_j}=
\sigma^*\left(\chi_{k_j}f_{m;j-1,k_{j-1}}{\cal B}_{m;j-1,k_{j-1}}\right),\,\,\,\,\,
{\cal B}_{m;j,k_j}=\sigma^*(\chi_{m,k_j}{\cal B}_{m;j-1,k_{j-1},l})
\end{eqnarray*}
and 
\[
 I_{\ell;m;\ell}[P](f)=\sum_{k_\ell=1}^{\nu_m}\int_{N_{\ell,k_\ell}}
\left({\cal H}^2(P_{m;\ell,k_\ell})\right)^{\ell+1}\left(f_{m;\ell,k_\ell}
{\cal B}_{\ell,k_\ell}\right)^2  
\]
Then we derive the recursive formula for these integrals.
and then through the inequality, which is derived thorugh Young's inequality
and for \( |x_{k_\ell}|<1\)
\[
 \left({\cal H}^2(P_{m;\ell-1,k_{\ell-1}})\right)^{\ell}\leq 
C_\ell m^{2\ell}\left(\frac{1}{x_{k_\ell}^{4\ell(\ell+1)}}+
\left({\cal H}^2(P_{m;\ell,k_\ell})\right)^{\ell+1}\right)
\]
then 
\begin{eqnarray*}
 I_{\ell-1;\ell-1}[P](f)\leq C_{\ell,m,\epsilon}\left(
||\nabla f||_{H^{2\ell(\ell+1)}(N_m)}
+ I_{\ell;\ell}[P](f)\right)
\end{eqnarray*}
The first term originates from the (EGHI), Leibniz's rule and  
the application of Hardy's inequality for the term 
\(D_{k_j}{\cal B}_{m;j,k_j}\) to transfer the derivative to the \(f\)-term.
The blow down process in the conical charts gives us the Sobolev space norms.
Specifically we have that 
\[
 |D_{x_{k_\ell}}^jf_{m;\ell,k_\ell}|=|\left(\frac{1}{x_{k_\ell}}
E\right)^j(\sigma^*
\left(\chi_{k_j}f_{m;j-1,k_{j-1}}{\cal B}_{m;j-1,k_{j-1}}\right))|
\]
Therefore by expanding the differentiation we arrive at 
\begin{eqnarray*}
 |D_{x_{k_\ell}}^jf_{m;\ell,k_\ell}|^2\leq C
\sum_{i_1+i_2+i_3=j}\frac{1}{r^{2i_1}}|\nabla^{i_2}\sigma^*\left(
\chi_{k_j}f_{m;j-1,k_{j-1}}|^2|\nabla^{i_3}{\cal B}_{m;j-1,k_{j-1}}^2\right)\leq\\
\leq C'\sum_{i_1+i_2=j}\frac{1}{r^{2i_1}}|\nabla^{i_2}f_{m;j-1,k_{j-1}}|^2
{\cal B}_{m;j-1,k_{j-1}}^2
\end{eqnarray*}
We iterate this process and finally we end up with the following inequality
\begin{eqnarray*}
I_{m}[P](f)\leq C\left(||\nabla f||_{H^{\beta_m}(N_m)}+
I_{m;1}[P](f)\right)
\end{eqnarray*}
where the term
\[
 I_{m;1}[P](f)=\sum_{k_{\gamma}=1}^{\nu_m}\int_{N_{m;\gamma_m,k_{\gamma_m}}}
\left({\cal H}^2(P_{m;\gamma_m,k_{\gamma_m}})\right)^{\gamma_m+1}
(f_{m;\gamma_m,k_{\gamma_m}})^2
\]
contains the polynomial mith mulitplicity \(m-1\) after the first generation of
 \(\gamma_m\) years of blow-ups for the \(m\) stratum and 
\( \beta_m=2\sum_{\ell=1}^{\gamma(m)}(\ell^2+\ell)=2S_2(\nu_m)+S_1(\nu_m)
\sim \nu_m^3\).

\paragraph{Conclusion} Having exhausted  the first generation of 
\(\gamma_m\)-years of blow -ups we reduced the multiplicity to \(m-1\)
at the cost of bringing in the Hardy - factor in the \(\gamma_m+1\) power, and we 
integrate in the tubular neighbourhood of set of lower multiplicity. Again the 
''center defining maps'' enter and are composed and we proceed.
Higher order derivatives appear but they are treated due to the formula that 
we calculated above.

%% file: ghi_proof2.tex
\subsection{The inequality for homogeneous polynomials.}

Let $P\in {\cal P}^{gH}$ then the inequality receives the simple form with 
the first derivatives. We start blowing up the  
origin which supports the maximal multiplicity of $V$ we obtain 
inequalities corresponding to the traces of the algebraic set on the balls
\[
\mpalla{0,\alpha,n}= V(x_j-1)\cap \kw{\alpha;j,n+1},\,\,\,\,\,\,\,\,
 j=1,\dots,n+1
\] 
Then we use the preceding inequalities for inhomogeneous polynomials and by a
scaling transformation we obtain the desired result by trowing the ``higher 
order terms''.  In the sequel we assume without loss of generality that we 
are already reduced to  the essential variables and proceed through the 
conical partition of unity and the notation of the preceding paragraph: 
\[
I_d[P](f)=\int_{\xw{n+1}} {\cal H}^1(P)f^2(x) v_{n+1}= 
\sum_{j=1}^{\nu_d}\int_{\ankw{\alpha;j,n+1}} 
(f_{d;1,j})^2x_j^{\nu_d-3}(P_{d;1,j})^{-\frac2d}\hat{v}_{n+1}
\]
Now in each cone \( \ankw{\alpha;j,1}\) the function 
\(f_{d,1,j}=\sigma^*(\chi_jf)\) is compactly supported and wehave that the volume 
form decomposes for the volume form of \( \xw{n}, v_n\)
\[
 \widehat{v}_{n+1}=dx_j v_n
\]
and the integral splits as 
\begin{eqnarray*}
 I_d[P](f)=\sum_{j=1}^{\nu_d}\int_{\xw{}} 
I_{d-1}[P_{d;1,j}](f_{d;1,j})(x_j)x_j^{\nu_d-3}dx_j 
\,\,\,\,\, (\mbox{CS})
\end{eqnarray*}
In each term we have \(\mbox{GHI}_1\) in the cross section: we obtain for
\(\tilde{\kl}_j\) the gradient in all but the \(x_j\)-variable:
\[
I_{d-1}[P_{d;1,j}](f_{d;1,j})\leq C  \int_{\mpalla{0,\alpha,n}}\sum_{i=0}^{h_{1,j}}
|\tilde{\kl}^i f_{d;1,j}|^2 = \left(R_{d,1,j}[P](f_{d;1,j})(x_j)\right)^2
\]
where the function is compactly supported in \(\xw{}\) and smooth with respect to 
\(x_j\). Hence we return back to (CS) and we have that 
\[
I_d[P](f)\leq C \sum_{j=1}^{\nu_d}
\int_{\xw{}}\left(R_{d,1,j}[P](f_{d;1,j})(x_j)\right)^2x_j^{\nu_d-3}dx_j
\]
We apply the  Hardy's inequality in each summand and obtain: 
\[
\int_{\xw{}}\left(R_{d,1,j}[P](f_{d;1,j})(x_j)\right)^2x_j^{\nu_d-3}dx_j\leq 
\frac{16}{(\nu_d-1)^2}
\int_{\xw{}}\left|\partial_{x_j}R_{d,1,j}[P](f_{d;1,j})(x_j)\right|^2
x_j^{\nu_d-1}dx_j
\]
Cauchy-Schwarz inequality suggests that 
\[
 \left|\partial_{x_j}\left(\sum_{i=1}^N G_i^2\right)^{1/2}\right|=
\left|\frac{1}{\sum_{i=1}^NG_i^2}\left(\sum_{i=1}^NG_i\partial_{x_j}G_i\right)\right|
\leq \left|\sum_{i=1}\partial_{x_j}G_i\right|
\]
Then apply this for \(G_i(x_j)=\left(\int_{\xw{n}}f^2(\cdot,x_j)\right)^{1/2}\) 
and thanks to Cauchy-Schwarz again we pass under the integral to obtain thanks to
Leibniz' srule that in gives us the following 
\[
 \int_{\xw{}}\left|\partial_{x_j}R_{d,1,j}[P](f_{d;1,j})(x_j)\right|^2
x_j^{\nu_d-1}dx_j\leq \int_{\ankw{\alpha;j,n+1}}\sum_{i=1}^{h_{1,j}}
|\kl^i f|^2 x_j^{\nu_d-1}\hat{v}_{n+1}
\]
Then we blow down in each conical chart and taking into account the usual 
estimates (we apply Hardy's inequlaity for the conical atlas partition of
unity) and  obtain that
\[
I_d[P]{f}\leq C 
\sum_{j=1}^{n+1}\int_{\ankw{\alpha;j,n+1}} 
\left(\sum_{i=1}^{h_{i,j}}|\kl_if|^2\right)v_{n+1}
\]
In the end we scale $x \mapsto {\tilde x}:= (\lambda^{-1} x_1,\dots,\lambda^{-1} x_{n+1})$
 Evidently, we let $\lambda\rightarrow \infty$ to obtain the inequality.
Similarly we obtain the inequality for \({\cal H}^2(P)\). q.e.d.

\noindent
{\bf Remark.} Sticking on the inequality for 
${\cal H}^2=\left|\frac{\kl P}{P}\right|^2$
one could also deduce the following inequality, actually 
by partial integration:  for $C_3(P)$ we have 

\[ 
 \olwn\left|\frac{\Delta P}{P}\right|f^2  \leq C_{3}(P) \olwn|\kl f|^2
\]

\subsection{ Further inequalities.}
Let \( P\in{\cal P}^{gH}\) of degree \(d\) and for  
\( s\in \xw{}\) then the following more general 
inequality is true for \( f \in C_0^\infty(\xw{n}\setminus V(P)) \):
\[
\olwn |P|^{s} f^2 \leq C(P,s) \olwn |P|^{s+\frac2d}|\kl f|^2 
\]
First we observe that the homogeneity we localize near \(V(P)\), 
the rest being treated by rectilinearization and partial integration.
We have for  
$f \in C^{\infty}_0(\xw{n}\setminus V(P)), \beta= s+\frac{2}{d}$:
and Young's inequality for \(p,q>1,pq=p+q\):
\[
 |P|^sf^2\leq \frac{\varepsilon^p|P|^{p\beta}f^2}{p}+
\frac{1}{\varepsilon^qq}|P|^{-\frac{2q}{m}}f^2,\,\,\,\,\,\,\,\,\,\,
(\mbox{BI}_1)
\]
then we have that 
\begin{eqnarray*}
\olwn |P|^2f^2 \leq \olwn \frac{1}{p}|P|^{p\beta}f^2 +
\frac{|P|^{-\frac{2q}{m}}f^2}{q} \leq \\
\leq \olwn \frac{1}{p}|P|^{p\beta}f^2
+C\olwn \left|\kl \left(P^{-\frac{2(q-1)}{m}}f\right)\right|^2
\end{eqnarray*}
Then the last term gives us that for \(\kappa_1=\frac{(q-1)^2}{m^2} \)
\[
\olwn |\kl \left(P^{-\frac{2(q-1)}{m}}f\right)|^2\leq 2 \left(\kappa_1 
I_1+I_2\right)
\]
where 
\[
I_2=\olwn |P|^{-\frac{2(q-1)}{m}} \left|\kl f\right|^2
\]
Now  
\[
I_1= \olwn P^{-\frac{2(q-1)}{2m}}{\cal H}^2(P)f^2\leq C_1(P)\left[(1+\eps)
\kappa_1I_1+(1+\frac1\eps)I_2\right] 
\]
Now we selct \(q\) such that:
\[
C_1\kappa_1(1+\eps)\leq \frac{1}{1+\eps}\Rightarrow q\leq 
1+\frac{m}{C_1(1+\eps)}
\]
and we have that
\[
I_1\leq C_3 I_2
\]
and finally we have that:
\[
\olwn |P|^sf^2 \leq C\left(\olwn |P|^{2\beta p}f^2+\olwn
|P|^{-\frac{2(q-1)}{m}} |\kl f|^2 \right)
\]
Now we split  \(N_\eta(P)=\{ x\in\xw{n}/ |P(x)|\leq \eta^m\},  \eta<1 \) 
in sets 
\[
N_i=\{ x\in \xw{n}/  \eta^{m(i+1)} \leq |P(x)|\leq \eta^{mi}\}
\]
and accordingly \( F_i=\fo(f)\cap N_i\):
\[
\olwn |P|^s f^2 =\sum_{i=1}^\infty \int_{F_i} |P|^sf^2
+\int_{\xw{n}\setminus N_\eta} |P|^sf^2 
\]
The second term is reduced to the 1-d problem by change of variable. The 
integral near the zero set is treated using the inequality derived 
from \((\mbox{BI}_1)\) for \(\varepsilon=\eta^{m(i+1}\) or
\(\varepsilon=\eta^{mi}\) analogously for \(\beta>0,<0\). Finally we scale 
and obtain tyhe result. q.e.d.

The inequality for the homogeneous polynomial allows to improve the
inequality for inhomogeneous polynomials, keeping only first derivatives:

\begin{inhom}
 Let $P : \xw{n} \rightarrow \xw{}$ be a polynomial
function of degree $d$ from the class ${\cal P}^{gH}$ with 
\(V(|\kl P|^2)\subset V(P)\). Let 
$f \in C^{\infty}_0({\bf R} ^n \setminus V(P))$  then for 
$C_3(P)$ it is true that:
\[
\int_{{\bf R} ^n} P^{-\frac2d}f^2 v_n \leq 
C_3(P)\int_{{\bf R} ^n} f^2 + (1 + |x|^2) |\kl f|^2 v_n 
\]
\end{inhom}

Proof. From the assumption we localize near \(V(P)\), the rest being treated
by rectilinearization and partial integration. 
We stratify the  algebraic set $V(P)$ by multiplicity:
\[
 V(P)=\Sigma_1\bigcup \dots \bigcup \Sigma_\mu
\]
Since \(P\) is inhomogeneous then $d > \mu$. 
The smooth part being clear, we will study the singular strata.
Let \(a \in \Sigma_m, x(a)=0\) and \( P_{0,m}=\mbox{in}_a(P)\) 
and introduce the   tubular neighbourhood of 
$\Sigma_m$ of width $\eps^{\frac dm}$
\[
N_\eps=\{x\in \xw{n}/ P^2_{0,m}(x)+P^2(x)<\eps^{2d}\}
\]
We write this as  
\[
N_\eps= C_\eps\bigcup S_\eps \bigcup R_\eps
\] 
where for suitable \(\delta<1\):
\begin{eqnarray*}
 \ced=\{ x \in N_\eps :|P| \geq  \delta^2|P_{0,m}| \} \\ 
 \sed=\{ x \in N_\eps :|P| \leq  \delta |P_{0,m}| \} \\
 \red=\{ x \in N_\eps : \delta^2|P_{0,m}|\leq |P|<\delta |P_{0,m}|\}
\end{eqnarray*}
Notice that the stratum of multiplicity 
\( \Sigma_m\subset \sed, \ced, \red \). 

The  integral 
\[
I_\eps=\int_{N_\eps} P^{-2/d} f^2 
\]
splits in three parts , \(i=1,2,3\):
\[
I_i=\int_{N_\eps} P^{-2/d} \psi_i^2f^2 
\]
for the functions localizing in the sets \(\ced,\sed,\red\) respectively :
\begin{eqnarray*}
\psi_1^2=\varphi\left(\frac{\delta^2|P_{0,m}|}{|P|}\right)\\
\psi_2^2=\varphi\left(\frac{|P|}{\delta|P_{0,m}|}\right)\\
\psi_1^2=1-\phi_1^2-\phi_2^2 
\end{eqnarray*}
We would  estimate them seperately:

\begin{enumerate}

\item Estimate for \(I_3\). The crucial observation is that 
\( \fo(\psi_3f)\cap V(P_{0,m})=\emptyset \) hence we apply the inequality 
for \(P_{0,m}\):
\[
I_3 \leq \delta^{-2/d}\int_{\red} |\kl (\phi_3 f)|^2  
\]
Now in \(\red\) we have that for \(\kappa_1 <\frac12\) 
\[
 |\kl \psi_3|\leq  C_1 + \kappa_1 
\left|\frac{\kl P_{0,m}}{P_{0,m}}\right|^2 
\]
and we apply the inequality for \({\cal H}^2(P_{0,m})\) and we are done.

\item The integral \( I_2 \) is filtered in the sets, for \(j>1\):

\[
\sedj=\{ x\in N_\eps / \delta^{j+1}|P_{0,m}|< |P|< \delta^j|P|\}
\]
and we obtain through the cut-off estimate :
\[
|\kl \psi_2|\leq C+\left(\frac{\delta}{2}\right)^j
\left|\frac{\kl P_{0.m}}{P_{0,m}}\right|^2 
\]
and the fact that \( \fo(f)\cap V(P_{0,m}) \) therefore we have that 
\[
I_2 \leq  C\olwn f^2+ |\kl f|^2 
\]
\item The last is the term near the cone \(V(P)\). Then introduce the sets
for \(\eta<1\):
\[ 
\cehj= \{x\in N_\eps/ |P|\geq  \delta \eta^j \geq \delta |P_{0,m}|  \}
\] 
Then through the cut-off estamates that allow us to  select appropriately
\(\eps\),  
we conclude that:
\[
I_1\leq  \olwn f^2+|\kl f|^2
\]
\end{enumerate}

Summing up we find 
\[
I_{1,\eps}\leq C\olwn |\kl f|^2 +f^2  
\]

q.e.d

\subsection{Example: cubic hypersurfaces.}
Here, we treat the inequality for cubic forms $(n+1)$- variables, defining 
a conic variety with singular part containing a line through the origin. 
The case $n=2$ is treated in [P] with elementary methods. 

If we choose the $x_1-$ axis such that $V(P)$ is singular along it then 
the form is written as $P(x)=x_1Q(\tilde{x})+C(\tilde{x})$ where 
$\tilde{x}=(x_2,\dots,x_{n+1})$ and $Q,C$ are  quadratic and cubic forms 
without common factors since otherwise the inequality is fairly easier. 
The form $P$ defines a cubic hypersurface in $\pxw{n}$ with 
a singularity at the point  $E_1$ with homogeneous coordiantes 
\[
E_1\equiv [1:0:\dots:0]
\] 
as well as at the points of the variety $T$ where  the varieties $V(C),V(Q)$ 
touch each other or the singular set $\Sigma$ of $V(C)$ on the hypeplane
$H_1$. We assume also without loss of generality that $\mbox{rank}(Q)=n$.

For the sake of simplicity we treat first the smooth part of the hypersurface
defined for $\epsilon >0$ by 
\[
V_{\epsilon}(P)= \bigcup_j V_{\epsilon,j}(P), 
\,\,\,\,\,\, V_{\epsilon,j}= \{ x\in \xw{n+1}/ 
\left|\frac{\partial P}{\partial x_j}\right| > \epsilon r^2 \}
\] 
Observe that in 
$V_{\epsilon}(P)$ it is true that for $c_i>0$, depending on the values of   
$P$ on $S^n$ that 
$c_1r^2\leq |\kl P| \leq c_2 r^2, c_1=(n+1)\epsilon $ and also 
that 
$c_3r\leq \left|\frac{\partial^2 P}{\partial x_i\partial_j}\right| \leq c_4 r $. 
Therefore the integral $I_{1,\epsilon}(P)[f]=I(P)[\chi_{1,\epsilon} f]$ 
is studied in $V_\epsilon(P)$ by the change of variables for  each $j$ by 
\begin{eqnarray*}
\Psi_j: V_{\epsilon,j} \rightarrow {\bf R}^{n+1} ,\\
(x_1,\dots,x_{n+1})\mapsto (x_1,\dots,x_{j_1},\xi,x_{j+1}, \dots,x_{n+1}),
\xi=P(x)
\end{eqnarray*}
Then we calculate $v_{j,n+1}=dx_1 \dots dx_{j-1}d\xi dx_{j+1}
\dots dx_{n+1}$
\[
\int_{V_{\epsilon}(P)}P^{-\frac23}f^2 v_{n+1} =
\sum_j \int_{\Psi_j(V_{\epsilon,j})} \xi^{-\frac23}\tilde{f}^2 
\left|\frac{\partial P}{\partial x_j}\right|^2 v_{j,n+1}\leq 
C\int_{V_{\epsilon}} |\kl f|^2 v_{n+1}
\] 
 by applying the 
preceding inequalities in combination with that of Hardy. 

To proceed to $I_{2,\epsilon}(P)[f]$ further we blowup 
$\sigma_{n+1}:\widehat{\xw{}}^{n+1}\rightarrow \xw{n+1}$ 
through the usual cones for \(\alpha>1\):
\[
\kw{\alpha;j,n+1}=\{ x \in \xw{n+1}/ x_j^2 \geq \frac{1}{1+\alpha}|x|^2\}
\]
The procedure will incorporate an induction with respect to $n$. 
Explicitly we have the formulas:

Let $j\neq 1$ then after restriction to $\kw{\alpha;j,n+1}$, 
$\sigma_{n+1}$ maps as 
\begin{eqnarray*}
i\neq j: x_i=v_iv_j,  \\
x_j=v_j, \\
P(\sigma_{n+1}(v))=x_j^3.P_j^1,\\ 
P_j^1=x_1Q_j^1(\tilde{v})+C_j^1(\tilde{v})
\end{eqnarray*} 
where $\tilde{v}=(v_2,\dots,v_{j-1},v_{j+1},\dots,v_{n+1})$.  
At the point  $a\in T_j$ or $a\in\Sigma_j$ the order 
$\mbox{ord}_a(P_j^1)=2$ and we 
choose coordinates such that $u(a)=0$ and setting $\tilde{u}=(u_2,\dots,u_n)$ 
we assume that $P_j^1(u_1,0,\dots,0)\sim u_1^2$ hence
\[
P_j^1(u)=(1+l_j(u))u_1^2+b_j(\tilde{u})u_1+ R_j(\tilde{u})
\]
where 
\begin{eqnarray*}
\mbox{deg}(l_j)=1,\,\,\,\, \mbox{deg}(\mu_j)=1,\,\,\,\,
\mbox{ord}_a(l_j)=0,\mbox{ord}_a(b_j)=1
\end{eqnarray*}
as well as that $R_j$ is of the smae form as $P_j^1$ but in $(n-1)-$variables. 
Then we determine  inductively the center 
\[
K_{n,j}=\{x \in \mpalla{0,\alpha,n}/ \mbox{ord}_x(P_j^1)=2\}
\]
of blowups and the years \(\nu_2\) of blow-ups for the generation \(m=2\). 
The formulas for the blowups in the $j_1$-conical chart, $j_1\neq 1$ are:  
\[
u_j=v_jv_{j_1}, j\neq j_1, v_{j_1}=u_{j_1}
\] 
implying that 
\[
P_j^1(u)=v_{j_1}^2P_{jj_1}^2(v) , \hspace{0.2truecm}
P_{jj_1}^2(v)=(1+v_{j_1}l_{jj_1}(v))v_1^2+b_{jj_1}(\tilde{v}))v_1 + 
R_j(\tilde{v})
\]
These blowups are repeated up to the desingularization of 
$V(R_j)$ and at each step we apply Hardy's inequality for the factor coming 
out of the blowup; finally we obtain a smooth chart as precedingly and we 
apply Hardy's inequality.  The $j=1$ conical chart is carried through 
analogously.

In conclusion we have that 
\[
\int_{\xw{n+1}}P^{-\frac23}f^2 v_{n+1}\leq C
\int_{\xw{n+1}}|\kl f|^2v_{n+1}+
\sum_{\ell=2}^{\nu_2}\int_{\xw{n+1}}\left|\kl^\ell f\right|^2 v_{n+1}
\] 
and  scaling we get the desired result.

%% file: ghi_appl_1.tex
\section{Heat expansion for operators
$H_{c,\alpha}=-\Delta + c |P|^{-\alpha}$}

In this section we will establish the existence of the small time expansion
 of the heat trace of the operator 
\[
H_{c,\alpha}=-\Delta + c |P|^{-\alpha}
\]
Actually we will prove the \(\lambda\rightarrow \infty\) expansion 
of the distributional trace  \(\mbox{tr}(R_{c,\alpha}(\lambda)^k\chi)\)
of the \(k\)-th power of  
\[
R_{c,\alpha}(\lambda)=(H_{c,\alpha}-\lambda)^{-1}
\]
This will achieved in the following steps:
\begin{enumerate}

\item Determination of the domain of self-adjointness of \(H_{c,\alpha}\)

\item Estimates  for 
\(R_{c,\alpha}(\lambda)=(H_{c,\alpha}-\lambda)^{-1}\) in   various operator 
norms

\item The preceding allow us to prove the expansion by considering the
expansion of each term of the Neumann series for \(R_{c,\alpha}(\lambda)\) 
around \(R_0(\lambda)=(-\Delta_n-\lambda)^{-1}\). 

\item After identifying the form of these terme we appeal to the usual 
Mellin transform theorem in view of the Atiyah-Bernshtein-Gelfand theorem 
on the meromorphic continuation of integrals depending on complex powers

\end{enumerate} 

\paragraph{Domain of Self-Adjointness}  
First we determine the domain of selfadjointness of  the operator  
$H_{c,\alpha}=-\Delta + c |P|^{-\alpha}$ where 
$c \in C_{0}^{\infty}({\bf R}^n)$ and $P$ is a homogeneous polynomial of 
degree $d$ in the class \( {\cal P}^{gH}\). 
We extend $C_{0}^{\infty}({\bf R}^{n})$ provided that
the exponent $\alpha < \frac nd$.

The Kato-Rellich theorem (\cite{rs}) suggests   
that  the operator  $H_{c,\alpha}$ is essentially self adjoint on 
$C_{0}^{\infty}(\xw{n})$. The neccessary semiboundedness estimate
reads as 

\begin{essesa}
  Let $c \geq 0, \phi \in C_{0}^{\infty}(\xw{n}), 
\alpha<\frac nd, \kappa <1$ then it is true that 
\[\parallel c. |P|^{-\alpha} \phi \parallel _{L^2} \leq \kappa 
\parallel -\Delta \phi \parallel _{L^2} + 
\beta  \parallel \phi \parallel_{L^2}\]
Furthermore the domain closure of the operator 
$H_{c,\alpha}|(C_{0}^{\infty}({\bf R}^n))$ for small $\alpha$
consists of  $H^2({\bf R}^n)$.

\end{essesa}

Proof. We decompose the integral in two parts, through a suitable partition 
of unity: one sufficiently close - at distance \(\epsilon^{1/d}\)-
to the algebraic set and one in the complement.  Set then 
\[
N_\epsilon= \{ x \in \xw{n} / |P(x)| < \epsilon \}
\]
and 
\[
K_\epsilon = \xw{n}\setminus N_\epsilon
\]
Choose $\chi_{1,\epsilon}^2 + \chi_{2, \epsilon}^2=1$,
$\fo(\chi_{1,\epsilon}) \subset K_{\epsilon}$ the functions that execute 
this decomposition (it just suffice to take functions of the form 
\(\chi\left(\frac{P}{\eps}\right)\)) Split then 
\[
\phi^2 = \phi_{1}^2 + \phi_{2}^2=
(\chi_{1,\epsilon}\phi)^2 + (\chi_{2,\epsilon}\phi)^2
\]
 Then for the part localized in $K_{2, \epsilon}$, we use the inequality :
\[
|P|^{-\alpha} \leq \varepsilon P^{-2/d} + B_{\varepsilon}
\]
Hence, we compute and obtain
\begin{eqnarray*}
\parallel c |P|^{-\alpha} \phi \parallel _{L^{2}} =
\parallel c.|P|^{-\alpha}{ \phi_{1}} \parallel _{L^{2}}+
\parallel c. |P|^{-\alpha} \phi_{2} \parallel _{L^{2}} \\
\leq M.\parallel \phi \parallel_{L^{\infty}}
\parallel \chi_{1,\epsilon} \parallel_{L ^{2}} +
\nu \parallel (-\Delta) \phi_{2} \parallel_{L ^{2}} +
\beta \parallel \phi \parallel_{L^2}
\end{eqnarray*}
having applied twice the inequality $P^{-2/d} \leq C (-\Delta)$
for the second term where the first term is due to the fact for
\(\alpha<\frac nd\) then \(|P|^{-\alpha}\in L^1_{\mbox{loc}}(\xw{n})\).
 Finally, we let $\epsilon \rightarrow 0$ and obtain the
desired inequality. The Kato - Rellich theorem gives the essential self-adjointness
we have been looking for; the operator $H_{c,\alpha}$ is bounded
below by $-\beta$. q.e.d

Further we will examine the closure of the operator, initially considered
on $C_{0}^{\infty}({\bf R}^n)$. For this we write :
 \begin{eqnarray*}
\parallel H_{c,\alpha} \phi \parallel_{L^2}^{2}=
\parallel -\Delta \phi  \parallel_{L^2}^2 +
\parallel c|P|^{-\alpha} \phi \parallel_{L^{2}}+\\
+(\phi,[(-\Delta) \cdot c|P|^{-\alpha}|P|^{-\alpha}+ c|P|^{-\alpha}\cdot
(-\Delta)]\cdot\phi) 
\end{eqnarray*}
Therefore combining the obvious inequality for the inner product with 
the inequality obtained above, we get the desired estimate.

\paragraph{ Operator estimates.}
For the Neumann series we will need the following estimates relative to the 
resolvent $R_{0}(\lambda)=((-\Delta)-\lambda)^{-1}$:

\begin{opes}
Let $\lambda $ be sufficiently large and outside a 
cone shaped region enclosing the positive real axis then the operator norm
is,
\begin{eqnarray*}
|| |P|^{-\alpha}R_{0}(\lambda) ||_{L^2}=0(|\lambda|^{-1 +
\frac{d\alpha}{2}}),\\
|| |P|^{-\alpha}R_{0}(\lambda) \partial_{i}||_{L^2}= 
0(|\lambda|^{-\frac12 + \frac{d\alpha}{2}})\\
||[|P|^{-\alpha}, R_0(\lambda)]||_{L^2}=
0(|\lambda|^{-1 + \frac{d\alpha}{2}}),\\
||[...[|P|^{-\alpha}, R_0(\lambda) ],...],R_0(\lambda)]||_{L^2}= 
0(|\lambda|^{-k + \frac{d\alpha}{2}})
\end{eqnarray*}

\end{opes}

Proof. Compute 
\[
(|P|^{-\alpha}R_{0}(\lambda)\phi, \,\,\,\,|P|^{-\alpha}R_{0}(\lambda)\phi)
\] 
with $\phi \in L^2(\xw{n})$.  
For this use a suitable partition of unity as in Proposition (2),
splitting  $\psi = R_{0}(\lambda)\phi \in C^{\infty}_{0}({\bf R}^{n})$
as before, far and close to the algebraic set $V(P)$. 
For the term involving  $\psi_{1}= \chi_{1} \psi$ apply Cauchy-Schwarz,
while for $\psi_{2}= \chi_{2} \psi$ use the inequality as follows:  
\begin{eqnarray*}
 ab \leq \frac {a^{p}}{p} + \frac {b^{q}}{q}, \hspace{0.5truecm} 
\frac1p+\frac1q=1,\\
 p=\frac{2}{d\alpha}, a= \delta^{\frac1p}P^{-2\alpha}, b= \delta^{-\frac1p},
\end{eqnarray*}
and for the two terms we have that 
\[
\parallel P^{-\alpha} \psi_{2}\parallel_{L^2} \leq
C_{\epsilon} \parallel \phi_{2} \parallel 
\]
 and also that
\[
\parallel \psi_{2} \parallel_{L^2} \leq C_{\epsilon}|\lambda|^{-\frac12}
 \parallel \phi_{2} \parallel 
\] 
Choosing then
\[
\delta = \delta ^{-\frac{\alpha d}{2-\alpha d}}|\lambda|^{-1}
\]
we get the desired estimate.The other estimates follow
in the same fashion as the first one. q.e.d.

\paragraph{The Neumann series}
We will prove the existence of an asymptotic expansion for a suitable power
 \[
R_{c,\alpha}^{k}(\lambda)=(\lambda - H_{c,\alpha})^{-k}
\]
and proving the existence for each term there.
Actually, $R_{c,\alpha}^{k}(\lambda)$ is represented by an integral kernel
since as we intend to show, it is in the Hilbert-Schmidt ideal.

The resolvent expansion for $R_{c,\alpha}(\lambda)$ is
\[
R_{c,\alpha}(\lambda) = \sum_{j=0}^{N-1}(R_{0}(\lambda)|P|^{-\alpha})^{j}.
R_{0}(\lambda)+O(|\lambda|^{-(1-\frac{d\alpha}{2})N-1})
\]
the remainder is considered in the $L^{2}$-norm.
Use the fact that 
\[
R_{c,\alpha}^{k}(\lambda)=\frac {1}{k!}
\pal{k}(R_{c,\alpha}(\lambda))
\]
to calculate the form of the Neumann series for \(R^k_{c,\alpha}(\lambda)\).
We set then  
\[
\Pi^{j}(\lambda)= (\resoul{}|P|^{-\alpha})^{j},
\]
so that 
\[
R_{c,\alpha}(\lambda)= \sum_{i=0}^{N-1}\Pi^i \resoul{} +
+0(|\lambda|^{-(1-\frac{d\alpha}{2})N-1}) 
\]
In order to obtain the Neumann series of \(R_{c,\alpha}(\lambda)^k\)
to study the beaviour of the terms in this series.
Therefore we set us 
\[
I^{i,j}=\pal{j}\left(\Pi^{i}R_{0}(\lambda)\right)=
\sum_{\ell=0}^j (-)^{j-\ell+1}
\left( \pal{\ell}\Pi^i\right) \resoul{j-\ell+1}
\]
Therefore we need 
\[
\pal{\ell}\Pi^i=\sum_{m_1+\cdots+m_{\ell+1}=i,m_1\geq 1}
\Pi^{m_1}\resoul{}\Pi^{m_2}\resoul{}\cdots\resoul{}\Pi^{m_\ell+1}
\]
Here observe that the number of \(R_0(\lambda)\)-factors is \( \ell\)
Conclusively we have that  for the operator product factor 
\[
C_{i;\ell;m_1,\dots,m_{\ell+1}}=
\Pi^{m_1}\resoul{}\Pi^{m_2}\resoul{}\cdots\resoul{}\Pi^{m_\ell+1}
\]
\[
I^{i,j}= (-1)^{j+1}
\sum_{0\leq \ell \leq j,m_1\geq 1,  m_1+\cdots+m_{\ell+1}=i} 
C_{i;\ell;m_1,\dots,m_{\ell+1}} 
\cdot \resoul{j-\ell+1} 
\]
Now we see that 
\begin{eqnarray*}
||C_{i;\ell;m_1,\dots,m_{\ell+1}} ||_{L^2}\leq
|\lambda|^{-(1-\frac{d\alpha}{2})i-\ell}\\
||I^{i,j}||_{L^2} \leq |\lambda|^{-(1-\frac{d\alpha}{2})i-j-1} 
\end{eqnarray*}
The trace class norm estimates are derived through
the following estimate from \cite{c1} :

\begin{opesII}
 
Suppose $k>\frac{n(ql-l+2)}{4lq},l \geq 2, \chi \in 
L^{2p}_{0}(\xw{n}),\frac1p+\frac1q=1$. Then $\chi.\resoul{}$
is in $C_{l}(\xw{n})$ for $\lambda \in {\bf C}\setminus{\bar {\bf R}_+}$ and
\[
|| \chi \resoul{k}||_{l} \leq C_{l,\epsilon}
||\chi|| _{2p} |\lambda|^{-k + \frac{n(ql-l+2)}{4lq}}
\]
for $\fo(\chi) \subset K  \subset\subset {\bf R}^n$ for all $\lambda$ in 
$|Im \lambda| > \epsilon Re \lambda + \epsilon$ for given $\epsilon >0$.

\end{opesII}

Combining the preceding estimates (prop (9),(10)) we 
can conclude that $I^{i,j} \in C_{l}(\xw{n})$ and furthermore that 
\[
\parallel \chi I^{i,j}(\lambda)\parallel_{l} \leq C_{l,\epsilon,j}
\parallel \chi \parallel _{2p} |\lambda|^{-j + \frac{n(ql-l+2)}{4lq} -
i(1-\frac{d\alpha}{2})}
\]

\paragraph{The existence of the Asymptotic Expansion for 
$\mbox{tr}(\chi R_{c,\alpha}^{k}(\lambda))$.} 
We conclude this section  with the  proof of  the existence of the 
$\lambda \rightarrow \infty$ asymptotic expansion of the trace 
$tr(\chi \cdot R^{k}_{c,\alpha}(\lambda))$.
The latter for $k$ chosen conveniently is a sum of integrals of the form 
\begin{eqnarray*}
 {\cal T}^{j,k}(\lambda,\alpha_1,\dots,\alpha_j):=
\mbox{tr}(\chi{\cal I}^{j,k})(\lambda,\alpha_{1},...,\alpha_j)=\\
=\int_{{\bf R}^{nj}} \Pi_{i=1}^{j}
(R_{0}^{k_i}(\lambda)(x_i,x_{i+1}))
\cdot\Pi_{i=1}^j |P|^{-\alpha_i}(x_i)\chi(x)\omega_1  \cdots\omega_j
\end{eqnarray*}
with the convention $x_{j+1}\equiv x_1, 
X=(x_1,\dots,x_j) \in {{\bf R}}^{nj}$. The resolvent is given by 
\[
R_0(\lambda)(x_1,x_2)=
\frac{(\sqrt{-\lambda})^{\frac n2+2}}{2^{\frac n2}
|x_1-x_2|^{\frac n2 -1}}K_{n-\frac12}(\sqrt{-\lambda}|x_1-x_2|)=
\frac{1}{(4\pi)^{\frac n2}}\int_{0}^{\infty}
e^{\lambda y-\frac{(x_1-x_2)^2}{4y}}
\frac{dy}{y^{\frac n2}}
\]
The existence of the asymptotic expansion of the trace  
${\cal T}^{j,k}(\lambda,\alpha_1,\dots,\alpha_{j})$ is proved using the 
Mellin transform suggested by the following classical result - stated in
the terminology of \cite{c2}:

\begin{salmv}

 Let $f \in C_0^{\infty}({\bf R}_+), m\in {\bf R}$ such that 
$x^m f \in L^1_{loc}({\bf R}_+)$. Suppose that the Mellin transform 
$\widehat{f}$ has a meromorphic extension from the half plane 
$\{ z \in {\bf C} / \Re z > - m\}$ with poles and multiplicity function 
$S, d=deg(S)$ and further that 
$lim_{\Im s \rightarrow \infty}(|s|^{d}\widehat{f}(s))=0.$ Then 
$f \in \Gamma^{\infty}(\overline{{\bf R}_+})$.

\end{salmv}

Using the integral representations for the resolvent and  changing variables we
deduce that the Mellin transform of ${\cal T}_{j,k}$ with respect to 
$\lambda$ takes the form for $k=k_1+\cdots +k_j$ :
$$\widehat{{\cal T}}^{j,k}(s,\alpha_1,\dots,\alpha_{j})=c_{n,j,k}(s)
\int_{ {\bf R}^{nj} \times \overline{{\bf R}_+}^j}
{\cal P}_{(s,\alpha_1,\dots,\alpha_{j})}(x,y)\chi(x,y)dxdy $$ 
where we have set above 
\[
c_{n,j,k}(s):=4^{s-\frac{(2k-n)j}{2}-j}
\Gamma(s)\Gamma(-s-\frac{(2k-n)j}{2}-j)
\]
$y:=(y_1,\dots,y_j)$ and the integrand consists of:
\[
{\cal P}_{(s,\alpha_1,\dots,\alpha_j)}(x,y)=
\frac{\prod_{i=1}^{j}P^{\alpha_i}(x_i)
(\prod_{i=2}^jy_i^{s-j(k_i-\frac n2)-1}
({\tilde{\cal Q}}(x,y))^{-s+(\frac n2-k)j+1}}
{(1+y_2+\cdots+y_j)^{-s}}
\]
and  ${\tilde {\cal Q}}=i^*{\cal Q}$ is the restriction  of the 
homogeneous polynomial function 
${\cal Q}:{{\bf R}}^{nj} \times {\bf R}^j \rightarrow {\bf R} $ 
on the $y_1=1$ hyperplane
\[
{\cal Q}(x,y)= \sum_{k=1}^j(y_k+y_{k-1})\pi_k(y) x_k^2
-\sum_{i=1}^{j} \pi_{(i,i+1)}^j(y)(x_{i-1}y_i + x_{i+1}y_{i-1})
\cdot x_i
\] 
where $\cdot$ denotes the Euclidean inner product and we
have set as well 
\[
\pi_i^j(y)=y_1 \dots y_{i-1}y_{i+1} \dots y_j ,  \,\,\,\,\,\,\,\,
\pi_{(i,i+1)}^j(y)=y_1 \dots y_{i-1}y_{i+2}\dots y_j 
\]
The existence of the asymptotic expansion is a direct consequence of the 
following theorem on the meromorphic continuation of integrals depending on 
complex parameters, \cite{bg}.

\begin{abg}

 Let $P_1,...,P_k$ be regular functions
on ${\bf R}^n$ and $\phi \in C_0^{\infty}({\bf R}^n)$, then the integral 
\[
I(\lambda_{1},...,\lambda_{k})=\int_{{\bf R}^n} |P_{1}|^{\lambda_{1}}...
|P_{k}|^{\lambda_{k}} \phi \omega_{n}
\]
can be continued as a meromorphic function
on the whole space of the complex variables
$\lambda_{1},...,\lambda_{k}$; at the same time the poles can 
be situated on a finite number of series of hyperplanes of the form
$a_{1}\lambda_{1}+...+a_{k}\lambda_{k} +b +s =0$ where $a_1,\dots,a_k,b$
are fixed nonnegative integers and $s$ runs through all the odd natural
numbers

\end{abg}

%% file: ghi_appl_2.tex
\section{Growth of inegrals in semialgebraic sets}

Let \(P\in {\cal P}^{gH}\) be a homogeneous polynomial and  consider the following 
semialgebraic sets:
\[
 N_P(\eta)=\{ x\in\xw{n}/ \eps_0\leq P(x)\leq \eta \}
\]
Then we consider the function \(\zeta , \fo(\zeta)\subset \mpalla{0,R}\)
that satisfies the gradient estimate for \(\gamma,\delta>0\):
\[
 |\kl \zeta|\leq \gamma |\zeta|+\delta
\]
If \(\varphi\in C_0^\infty(\xw{}_+)\) with \(\fo(\varphi)\subset 
[0,1+\varepsilon), \phi\equiv 1 \) in \([0,1]\)  then
\(\varphi\left(\frac{P}{\eta}\right)\) localizes in the space 
\( P<\eta \). We will follow the classical identities for the monotonicity 
formulas, cf. \cite{hl}.
Then we have the following elementary identity:
\[
 \eta\partial_\eta\phi=-\frac{P}{|\kl P|^2} \kl P\cdot\kl \varphi 
\]
we can differntiate the integral replacing \(\zeta\) by 
\(\varphi\left(\frac{\eps_0}{P}\right)\zeta\) in order to localize in 
\(N_P(\eta)\):
\[
 I_\zeta(\eta)=\olwn \phi\zeta^2 
\]
and get:
\[
 \eta\frac{d I_\zeta}{d\eta}=-\olwn \frac{P}{|\kl P|^2}\zeta^2\kl P\cdot \kl\varphi
\]
Then integrate by parts to obtain for \(Q=|\kl P|^2\):
\[
\zeta\frac{d I_\zeta}{d\eta}\leq \olwn 
\left(\frac{P^2}{Q^2}\left|\frac{\Delta P}{P}\right|+\delta +2\gamma\frac{P^2}{Q}
+\frac{|P|}{|\kl P|}\frac{|\kl Q|}{Q} \right)\varphi\zeta^2+\mbox{vol}(N_\eta)
\]
Then applying the conical Lojasiewicz inequality we obtain:
\[
 \frac{|P|}{|\kl P|}\leq c |P|^{1/m}\leq c\eta^{1/m}
\]
If we assume that \(Q\in {\cal P}^{gH}\) -it is not neccessary, since with partial 
integration we can  avoid the term involving \( \kl Q\) -
 then applying the generalized Hardy's 
inequality we arrive at the differential inequality 
\[
 \frac{dI_\zeta}{d\eta}\leq C\left(\gamma^2 \eta^{\frac2m-1}I_\zeta+\eta^{\frac1m}\right)
\]
Then we arrive at the conclusion:
\[
 I_\zeta(\eta)\leq c_1e^{c_2\gamma^2\eta^{2/m}} I_\zeta(\eta/2)
\]
Actually this estimate  combined with Harnack inequalities for semialgebraic sets
(proved in \cite{p2})  provide growth estimates for functions
on semialgebraic sets.